\documentclass[11pt]{amsart}
\usepackage[T2A]{fontenc}
\usepackage[utf8]{inputenc}
\usepackage[english]{babel}
\usepackage[margin=2.75cm]{geometry}

\usepackage{amsmath, amsfonts}
\usepackage{mathtools}
\providecommand\given{} 
\newcommand\SetSymbol[1][]{%
	\nonscript\:#1\vert
	\allowbreak \nonscript\:	\mathopen{}}
\DeclarePairedDelimiterX\Set[1]\{\}{%
	\renewcommand\given{\SetSymbol[\delimsize]}	#1} 

\usepackage{array} 
\usepackage{multirow}
\usepackage{rotating}
\usepackage{booktabs}
\usepackage{subcaption}

\usepackage{tikz}

\usepackage[pdffitwindow = true, pdfstartview = FitH,  pdftitle = {2-neighborly 0/1-polytopes of~dimension 7},
pdfauthor = {A.N. Maksimenko}
]{hyperref}

\usepackage{numprint} 
\newcommand{\nm}[1]{\numprint{#1}}

\usepackage[lined,ruled]{algorithm2e} 

\theoremstyle{plain}

\theoremstyle{definition}

\theoremstyle{remark}

\newcommand{\Z}{\mathbb{Z}}  
\DeclareMathOperator{\conv}{conv} 
\DeclareMathOperator{\cone}{cone} 
\DeclareMathOperator{\ext}{ext} 

\title{2-neighborly 0/1-polytopes of dimension 7}

\author{Aleksandr Maksimenko}
\thanks{Supported by the~project \textnumero\,1.5768.2017/9.10 of~P.\,G.~Demidov Yaroslavl State University within State Assignment for~Research.}
\address{Laboratory of~Discrete and Computational Geometry, P.G. Demidov Yaroslavl State University, ul. Sovetskaya 14, Yaroslavl 150000, Russia} 
\email{maximenko.a.n@gmail.com}

\usepackage{fancyhdr}
\begin{document}

\begin{abstract}
We give a~complete enumeration of~all 2-neighborly 0/1-polytopes of~dimension 7.
There are \numprint{13959358918} 
different 0/1-equivalence classes of~such polytopes.
They form \numprint{5850402014} combinatorial classes and \numprint{1274089} 
different $f$-vectors.
It enables us to list some of~their combinatorial properties.
In particular, we have found a~2-neighborly polytope with 14~vertices and 16 facets.
\end{abstract}

\maketitle

\section{Introduction}
A \emph{0/1-polytope} is a~convex polytope whose set of~vertices is subset of~$\{0,1\}^d$.
For~beautiful introduction to the world of~0/1-polytopes, we refer to Ziegler~\cite{Ziegler:2000}.
Some recent results can be found in \cite{Gillmann:2006}. 
The~classification of~all \numprint{1226525} different 0/1-equivalence classes of~0/1-polytopes of~dimension~5 was done by Aichholzer~\cite{Aichholzer:2000}.
Also he completed the classification of~6-dimensional 0/1-polytopes up to 12 vertices.
Recently, Chen and Guo~\cite{Chen:2014} computed the numbers of~0/1-equivalence classes
 of~6-dimensional polytopes for each number of~vertices $k \in [13, 64]$.
But~nowadays it is too hard to list all these $\approx 4.0 \cdot 10^{14}$ classes 
 and investigate their properties explicitly.
(If every polytope will ocupy 8 bytes, then all the database will take about 3 petabytes.)
Thus, it makes sense to focus on some interesting families of~0/1-polytopes.

A convex polytope $P$ is called \emph{2-neighborly} if any two vertices form a~1-face (i.e. edge) of~$P$.
There are at least two reasons for investigation of~2-neighborly 0/1-polytopes:
\begin{enumerate}
	\item Let $P_{d, n}$ is a~random $d$-dimensional 0/1-polytope with $n$ vertices. 
	In 2008, Bondarenko and Brodskiy~\cite{Brodski:2008} showed, that if $n = O(2^{d/6})$, 
	then the probability $\Pr(P_{d,n}\text{ is 2-neighborly})$ tends to 1 as $d \to \infty$.
	Similar results are obtained by Gillmann~\cite{Gillmann:2006}.
	
	\item Special 0/1-polytopes, such as the cut polytopes, the traveling salesman polytopes,
	 the~knapsack polytopes, the $k$-SAT polytopes, the 3-assign\-ment polytopes, the set covering polytopes
	 and many others have 2-neighborly faces with superpolynomial (in the dimension) number of~vertices 
	 \cite{Deza:1997, Maksimenko:2014, Maksimenko:2017}.
\end{enumerate}

We enumerated and classified all \numprint{13959358918} 
0/1-equivalence classes of~7-dimensional 2-neighborly 0/1-polytopes.
It enables us to investigate extremal properties of~these polytopes.
For example, we have found a~2-neighborly polytope with 14 vertices and 16 facets.
This is the~first known example of~a 2-neighborly polytope (except a~simplex)
 whose number of~facets is not greater than the~number of~vertices plus~2.

In~\cite{Aichholzer:2000}, Aichholzer stated the~question about the~maximal number $N_2(d)$ 
of~vertices of~a 2-neighborly $d$-dimensional 0/1-polytope.
He showed that $N_2(6) = 13$, $18 \le N_2(7) \le 24$, $N_2(8) \ge 25$, $N_2(9) \ge 33$, $N_2(10) \ge 44$.
We improve these estimations: $N_2(7) = 20$, $28 \le N_2(8) \le 34$, $N_2(9) \ge 38$, $N_2(10) \ge 52$.

The~entire database occupies about 1TB. The~part of~it (in particular, all 6-polytopes) and the~list of~all f-vectors are available at \url{https://github.com/maksimenko-a-n/2neighborly-01polytopes}.

\section{Enumeration of~2-neighborly 0/1-polytopes}  

Every 0/1-polytope is a~convex hull of~a set $X\subseteq\{0,1\}^d$.
Since the~natural way for defining a~0/1-polytope is the~defining its set of~vertices $X$,
in the~following we will frequently call $X$ by a~``polytope'',
having in mind the~convex hull $\conv(X)$.

We will use the~following trivial facts. Every 0/1-polytope $\conv(X)$, $X\subseteq\{0,1\}^d$, is the~convex hull of~a 0/1-polytope $\conv(X \setminus \{x\})$ and a~vector $x \in X$.
The~same is true for 2-neighborly polytopes. Let $P$ be a~2-neighborly polytope and $X = \ext(P)$ be its set of~vertices.
Then for every $x \in X$ the~polytope $\conv(X \setminus x)$ is also 2-neighborly.
Thus, we can enumerate 2-neighborly 0/1-polytopes iteratively, starting with a~polytope consisting of~a single point and adding one point each time.

Two polytopes are \emph{0/1-equivalent} if one can be transformed into the~other by a~symmetry
of~the~0/1-cube. More precisely, this transformation means the~using of~two operations: permuting of~coordinates and replacing some coordinates $x_i$ by $1 - x_i$ (\emph{switching}).
Thus, one 0/1-equivalence class can contain up to $2^d d!$ of~0/1-polytopes of~dimension $d$.
The~0/1-equivalence implies affine and combinatorial equivalences~\cite[Proposition 7]{Ziegler:2000}.

Every 0/1-vector $x = (x_1, x_2, \dots, x_d) \in \{0,1\}^d$ can be associated with the~binary number ``$x_1 x_2 \ldots x_d$''.
Thus, every 0/1-polytope $X \subseteq \{0,1\}^d$ can be naturally associated with the~sequence of~such binary numbers sorted in increasing order.
Let $\mathcal{C}$ be a~0/1-equivalence class (a set of~polytopes). A polytope $X \in \mathcal{C}$ is called a~\emph{representative} if it is lexicografically less than any other polytope $Y \in \mathcal{C}$.
\SetKwFunction{Represent}{representative}%
For a~given 0/1-polytope $Y$, the~appropriate \Represent{$Y$} can be found with a~straightforward branch and bound algorithm.

Therefore, for enumeration of~2-neighborly 0/1-polytopes we can iteratively use the~algorithm~\ref{alg:enumeration}.
In the~first step, $T_1$ contains only one polytope $\{(0,\dots,0)\}$.

\SetKwProg{Proc}{Procedure}{}{end} 
\SetKwProg{Fn}{Function}{}{end} 
\SetKwInOut{Input}{Input}
\SetKwInOut{Output}{Output}
\begin{algorithm}
	\caption{The~enumeration of~2-neighborly 0/1-polytopes} 
	\label{alg:enumeration}
	\SetKwData{In}{$T_n$} 
	\SetKwData{Out}{$T_{n+1}$} 
	\SetKwData{Dim}{d}
	\SetKwData{Vert}{n}
	\SetKwData{X}{X}
	\SetKwData{v}{v}
	\SetKwData{Y}{Y}
	\SetKwData{y}{y}
	\SetKwData{false}{false}
	\SetKwData{true}{true}
	\SetKwFunction{AddVertex}{enumerate\_2neighborly}
	\SetKwFunction{IsNeighb}{is\_2neighborly}
	\SetKwFunction{noEdge}{no\_edge\_0y}
	\Input{ the~dimension $d$, the~array \In of~2-neighborly 0/1-polytopes with $n$ vertices (every polytope is an array of~$n$ 0/1-vectors)}
	\Output{ the~array \Out of~2-neighborly 0/1-polytopes with $n+1$ vertices}
	\BlankLine
	\Fn{\AddVertex{\Dim, \In}}{
		\For{$X \in \In$}{
			\For{$v \in \{0,1\}^d \setminus X$}{
				\If{\IsNeighb{$X$, $v$}}{
					add \Represent{$X \cup \{v\}$} to \Out\;
				}
			}	
		}
		sort \Out and remove duplicates\;
		\Return{\Out}\;
	}
	\BlankLine
	\tcp{Let $\conv(\X)$ be 2-neighborly. Is $\conv(\X \cup \{v\})$ 2-neighborly?}
	\Fn{\IsNeighb{\X, \v}}{
		$Y \coloneqq \emptyset$\;
		\tcp{firstly, test edges for \v and $x \in X$}
		\lFor(\tcp*[f]{switch \v to 0}){$x \in \X$}{add $x \oplus \v$ to Y}
		\For{$y \in Y$}{
			\lIf{\noEdge{$Y$, $y$}}{\Return{\false}}
		}
		\tcp{test edges $\{x,y\} \subseteq X$}
		\For{$x \in \X$}{
			$w \coloneqq \v \oplus x$\;
			$Y \coloneqq \emptyset$\;
			\lFor(\tcp*[f]{switch $x$ to 0}){$y \in \X \setminus x$}{add $y \oplus x$ to Y}
			\For{$y \in Y$}{
				\If{$w \wedge y = w$}{
					\lIf{\noEdge{$Y \cup \{w\}$, $y$}}{\Return{\false}}
				}	
			}
		}
	
		\Return{\true}\;
	}
	\BlankLine
	\tcp{Isn't $\{0,\y\}$ an edge of~$\conv(\Y \cup \{0\})$?}
	\Fn{\noEdge{\Y, \y}}{
		$Z \coloneqq \emptyset$\;
		\For{$z \in \Y \setminus \y$}{
			\lIf{$z \wedge \y = z$}{add $z$ to $Z$}
		}
			\lIf{$\y \in \cone(Z)$}{\Return{\true}}	
		\Return{\false}\;
	}
\end{algorithm}

For testing 2-neighborliness of~a 0/1-polytope $X \subseteq \{0,1\}^d$ we use the~ideas described in~\cite[Sec. 2.2]{Aichholzer:2000}.
Let $v, w \in X$ and we want to check the~adjacency of~$v$ and $w$ in $\conv(X)$.
First of~all, we switch $X$ to $Y = \Set*{x \oplus v \given x \in X}$.
(Here, $\oplus$ is a~coordinatewise XOR operation.)
Hence, $v$ will be switched to $0$. 
It is easy to prove that the~vertices $0, y \in Y$ form an edge of~a polytope $Y$ iff they form an edge of~a polytope $Z = \Set*{z \in Y \given z \wedge y = z}$.
(Here, $\wedge$ is a~coordinatewise AND operation.)
Thus, we have to check whether $y$ is in the~conical hull
\[
\cone(Z) = \Set*{\sum_{z \in Z} \lambda_z z \given \lambda_z \ge 0}.
\]
Namely, vertices $0$ and $y$ form an edge in $Z$ iff $y \notin \cone(Z)$.
The~cheking of~$y \notin \cone(Z)$ can be done by solving the~corresponding linear programming problem.
We did it with COIN-OR Linear Programming Solver~\cite{COIN}.

We have run this algorithm on the~computer cluster of~Discrete and computational geometry laboratory of~Yaroslavl state university (\url{https://dcgcluster.accelcomp.org}). 
The~cluster has a~hundred 2.9GHz-cores. After several weeks of~computations we had got the~rezults collected in Table~\ref{tab:1b}. 
Every 0/1-vector $x \in \{0,1\}^d$, $d \le 8$, we store as a~1-byte integer. Thus, a~polytope with $n$ vertices occupies $n$ bytes and all the~database~--- about 173GB.

\begin{table}[tbh]
	\centering
    \begin{subfigure}[t]{0.4\textwidth}
	\caption{The~dimension 6}
	\label{tab:1a}
	\begin{tabular}{rr}
		\toprule
		\textbf{vertices} & \textbf{0/1-equivalence}\\
		& \textbf{classes}\\
		\midrule 
1          & \nm{1}\\
2          & \nm{6}\\
3         & \nm{16}\\
4         & \nm{94}\\
5        & \nm{445}\\
6       & \nm{2528}\\
7      & \nm{12359}\\
8      & \nm{47445}\\
9     & \nm{108220}\\
10     & \nm{110032}\\
11      & \nm{38221}\\
12       & \nm{3222}\\
13         & \nm{36}\\
		\midrule 
		total & \nm{322625}\\
		\bottomrule
	\end{tabular}
	\end{subfigure}
\hspace*{1cm}
    \begin{subfigure}[t]{0.4\textwidth}
	\caption{The~dimension 7}
	\label{tab:1b}
	\begin{tabular}{rr}
		\toprule
		\textbf{vertices} & \textbf{0/1-equivalence}\\
								& \textbf{classes}\\
		\midrule 
1	& \nm{1}\\
2	& \nm{7}\\
3	& \nm{23}\\
4	& \nm{191}\\
5	& \nm{1510}\\
6	& \nm{16373}\\
7	& \nm{183209}\\
8	& \nm{1985525}\\
9	& \nm{18565154}\\
10	& \nm{136197421}\\
11	& \nm{707274277}\\
12	& \nm{2345160234}\\
13	& \nm{4456209397}\\
14	& \nm{4284931624}\\
15	& \nm{1757834961}\\
16	& \nm{244831279}\\
17	& \nm{8967617}\\
18	& \nm{73512}\\
19	& \nm{180}\\
20	& \nm{3}\\
		\midrule 
		total & \nm{13962232498}\\
		\bottomrule
	\end{tabular}
\end{subfigure}
	\caption{The~number of~0/1-equivalence classes of~2-neighborly polytopes of~dimensions 6 and 7}
	\label{tab:1}
\end{table}

Our results for the~dimension 6 coincide with Aichholzer database~\cite{Aichholzer:2000}.
In addition, we enumerate all 2-neighborly 0/1-polytopes of~dimension 6 with 13 vertices.

\section{Evaluating of~combinatorial types and f-vectors}  

It is well known (see e.g.~\cite{Kaibel:2002}) that the~combinatorial type (face lattice) of~a~polytope $P$ with vertices $\{v_1, \dots, v_n\}$  and facets $\{f_1, \dots, f_k\}$ is uniquely determined by its facet-vertex incidence matrix $M = (m_{ij}) \in \{0,1\}^{k \times n}$, where $m_{ij} = 1$ if
facet $f_i$ contains vertex $v_j$, and $m_{ij} = 0$ otherwise.
Thus, polytopes are combinatorially equivalent iff their facet-vertex incidence matrices differ only by column and row permutations.

For every polytope in our database we computed its facet-vertex incidence matrix $M$ by using \emph{lrs}~\cite{lrs}.
This evaluation takes about 10 days on the~computer cluster with 32 cores.
After that, for every matrix $M$, we computed the~canonical form of~a vertex-facet digraph of~$M$ by using \emph{bliss}~\cite{bliss} (as it was done in~\cite{Firsching:2018}).
This evaluation takes about 2 days on the~computer cluster with 32 cores.
Having sorted canonical forms, we have splitted all the~polytopes into combinatorial equivalence classes.

For computing f-vector of~a polytope from its facet-vertex incidence matrix, we used  Kaibel\&Pfetsch algorithm~\cite{Kaibel:2002} and modified it for the~case, when the~number of~vertices is small (an~incidence matrix row can be stored in a~32-bit integer). The~computing of~f-vectors of~all polytopes took about two weeks on the~cluster.

\begin{table}[tbh]
	\centering
	\begin{tabular}{rrrr}
		\toprule
		\textbf{vertices} & \textbf{0/1-equivalence} & \textbf{combinatorial} & \textbf{f-vectors}\\
		& \textbf{classes} & \textbf{classes} & \\
		\midrule 
6 & \nm{237} & \nm{1} & \nm{1}\\
7 & \nm{334} & \nm{2} & \nm{2}\\
8 & \nm{102} & \nm{8} & \nm{5}\\
9 & \nm{10} & \nm{7} & \nm{4}\\
10 & \nm{1} & \nm{1} & \nm{1}\\
		\midrule 
		total & \nm{684} & \nm{19} & \nm{13}\\
		\bottomrule
	\end{tabular}
	\caption{Full-dimensional 2-neighborly 0/1-polytopes of~dimension 5}
	\label{tab:d5}
\end{table}

\begin{table}[tbh]
	\centering
	\begin{tabular}{rrrr}
		\toprule
		\textbf{vertices} & \textbf{0/1-equivalence} & \textbf{combinatorial} & \textbf{f-vectors}\\
& \textbf{classes} & \textbf{classes} & \\
		\midrule 
		7 & \nm{9892} & \nm{1} & \nm{1}\\
		8 & \nm{46813} & \nm{4} & \nm{4}\\
		9 & \nm{108178} & \nm{81} & \nm{32}\\
		10 & \nm{110029} & \nm{9651} & \nm{180}\\
		11 & \nm{38221} & \nm{17782} & \nm{411}\\
		12 & \nm{3222} & \nm{2730} & \nm{455}\\
		13 & \nm{36} & \nm{35} & \nm{34}\\
		\midrule 
		total & \nm{316391} & \nm{30284} & \nm{1117}\\
		\bottomrule
	\end{tabular}
	\caption{Full-dimensional 2-neighborly 0/1-polytopes of~dimension 6}
	\label{tab:d6}
\end{table}

\begin{table}[tbh]
	\centering
		\begin{tabular}{rrrr}
			\toprule
			\textbf{vertices} & \textbf{0/1-equivalence} & \textbf{combinatorial} & \textbf{f-vectors}\\
& \textbf{classes} & \textbf{classes} & \\
			\midrule 
8 & \nm{1456318} & \nm{1} & \nm{1}\\
9 & \nm{17588780} & \nm{6} & \nm{6}\\
10 & \nm{135330686} & \nm{419} & \nm{108}\\
11 & \nm{706996729} & \nm{4790131} & \nm{2090}\\
12 & \nm{2345138023} & \nm{271351237} & \nm{17113}\\
13 & \nm{4456209206} & \nm{1414858979} & \nm{66929}\\
14 & \nm{4284931624} & \nm{2487091476} & \nm{171289}\\
15 & \nm{1757834961} & \nm{1431813684} & \nm{303063}\\
16 & \nm{244831279} & \nm{231549854} & \nm{382319}\\
17 & \nm{8967617} & \nm{8872600} & \nm{282000}\\
18 & \nm{73512} & \nm{73444} & \nm{48988}\\
19 & \nm{180} & \nm{180} & \nm{180}\\
20 & \nm{3} & \nm{3} & \nm{3}\\
			\midrule 
total & \nm{13959358918} & \nm{5850402014} & \nm{1274089}\\
			\bottomrule
		\end{tabular}
	\caption{Full-dimensional 2-neighborly 0/1-polytopes of~dimension 7}
	\label{tab:d7}
\end{table}

The~results of~these computations are collected in Tables~\ref{tab:d5}--\ref{tab:d7}.
We enumerate only full-dimensional 0/1-polytopes, since any nonfull-dimensional 0/1-polytope is affinely equivalent to some full-dimensional one~\cite{Ziegler:2000}.

To give an idea of~the~magnitude of~the~obtained numbers, we give a~couple of~examples:
f-vector (13, 78, 266, 531, 603, 355, 84) consists of
\nm{2448144} combinatorial classes;
f-vector (9, 36, 82, 114, 97, 48, 12) consists of~one combinatorial class with \nm{5160979} 0/1-equivalence classes.


For every combinatorial type, we store its facet-vertex incidence matrix. If the~number of~vertices (columns of~the~matrix) is not greater than 16, one row of~the~matrix occupies 2 bytes.
The~average number of~facets (rows of the matrix) is 97. Thus, all combinatorial types occupy about 1TB.
The~part of~the~database (in particular, all 6-polytopes) is available at \url{https://github.com/maksimenko-a-n/2neighborly-01polytopes}.
The~full information can be requested from the~author (by e-mail).

\section{$d$-Polytopes with $d+3$ vertices}

Combinatorial types of~$d$-polytopes with $d+3$ or less vertices can be enumerated by using Gale diagrams~\cite[Chap.~6]{Grunbaum:2003}. For $d+1$ vertices there are only one $d$-polytope~--- a~simplex. For $d+2$ vertices, the~number of~combinatorial types is equal to the~number of~tuples $(m_0, \{m_1, m_{-1}\})$, where $m_0, m_1, m_{-1} \in \Z$, $m_0 \ge 0$, $m_1 \ge 2$, $m_{-1} \ge 2$, and $m_0 + m_1 + m_{-1} = d+2$~\cite[Sec.~6.3]{Grunbaum:2003}. The~appropriate polytope is 2-neighborly iff $m_1 \ge 3$, $m_{-1} \ge 3$. For small $d$, these tuples can be easily enumearted by hands.

The~combinatorial type of~every $d$-polytope with $d+3$ vertices is defined by the~appropriate reduced Gale diagram or wheel-sequence~\cite{Fusy:2006}. We don't list here the~properties of~these interesting objects, since it was done in~\cite{Grunbaum:2003,Fusy:2006,Maksimenko:2018}. The~results of~enumerating wheel-sequences by a~computer are collected in Table~\ref{tab:gale}. 
They coincide with the~first values of~the~sequence A114289: \url{https://oeis.org/search?q=A114289} and with the~Fukuda--Miyata--Moriyama collection of~$d$-polytopes for $d \le 6$~\cite{Fukuda:2013}.

\begin{table}[tbh]
	\centering
	\begin{tabular}{rrrrrrr}
		\toprule
		 & \multicolumn{3}{c}{\textbf{\boldmath $d+2$ vertices}} & \multicolumn{3}{c}{\textbf{\boldmath $d+3$ vertices}}\\
		 \cmidrule(lr){2-4} \cmidrule(lr){5-7}
		\textbf{\boldmath$d$} & all & 2-neighborly & 2-neighborly & all & 2-neighborly & 2-neighborly \\
		 &     & polytopes    & 0/1-polytopes &     & polytopes     & 0/1-polytopes\\
		\midrule 
		4 &  4 & 1 & 1 &   31 &   1 &   0 \\
		5 &  6 & 2 & 2 &  116 & 11 &  8 \\
		6 &  9 & 4 & 4 &  379 & 85 & 81 \\
		7 &12 & 6 & 6 & 1133 &423 & 419 \\
		\bottomrule
	\end{tabular}
	\caption{Combinatorial types of~$d$-polytopes with $d+2$ and $d+3$ vertices}
	\label{tab:gale}
\end{table}

For $d \le 7$, every combinatorial type of~a 2-neighborly $d$-polytope with $d+2$ vertices can be represented by a~0/1-polytope. Almost the~same is true for polytopes with $d+3$ vertices. The~exceptions are 4 polytopes and the~pyramids over them. 
The~first one is a~cyclic 4-polytope with 7 vertices. 
The~second can be represented by the~wheel-sequence (0, 1, 0, 1, 1, 0, 1, 0, 1, 1, 1, 1). It has f-vector (8, 28, 50, 44, 16) and its facet-vertex incidence matrix has two columns with 12 ones as opposed to other polytopes with the~same f-vector. 
The~third polytope can be represented by the~wheel-sequence (0, 1, 0, 1, 0, 1, 0, 1, 0, 1, 0, 1, 1, 1). It has f-vector (8, 28, 51, 47, 18) and its facet-vertex incidence matrix has a~column with 14 ones as opposed to other polytopes with the~same f-vector.
The~forth polytope is represented by the~sequence (0, 1, 0, 1, 0, 1, 0, 1, 0, 1, 1, 1, 1, 1). It has f-vector (9, 36, 80, 103, 72, 22) and its incidence matrix has no column with 12 ones as opposed to other polytopes with the~same f-vector.

\section{Polytopes with a~small number of~facets}

The~minimal and the~maximal numbers of~facets of~2-neighborly 0/1-polytopes listed in Table~\ref{tab:facets}.
As can be seen, there is a~2-neighborly 7-polytope $P_{14, 16}$ with 14 vertices and 16 facets. 
In Figure~\ref{fig:p14}, we list vertices of~$P_{14, 16}$.
As far as we know, any other 2-neighborly polytope (except a~simplex) has the~property $(\text{facets} - \text{vertices}) \ge 3$. The~polytope $P_{14, 16}$ has several other special properties. It is the~only 2-simple polytope in our database. (A $d$-polytope is \emph{2-simple} if every $(d-3)$-face is incident to exactly three facets.) All its vertex figures are combinatorially equivalent 6-polytopes with 13 vertices and 11 facets.
For any vertex figure of~any other polytope in our database, the~number of~facets is not less than the~number of~vertices.

\begin{table}[tbh]
	\centering
	\begin{tabular}{l*{13}r}
		 \multicolumn{6}{c}{dimension 5} \\
		\cmidrule{1-6}
vertices & 6 & 7 & 8 & 9 & 10 \\
facets min & 6 & 10 & 12 & 16 & 22 \\
facets max & 6 & 12 & 20 & 22 & 22 \\
		\cmidrule{1-6}
		 \multicolumn{8}{c}{dimension 6} \\
		\cmidrule{1-8}
vertices & 7 & 8 & 9 & 10 & 11 & 12 & 13\\
facets min & 7 & 11 & 13 & 14 & 17 & 21 & 26\\
facets max & 7 & 16 & 30 & 47 & 55 & 65 & 76 \\
		\cmidrule{1-8}
		 \multicolumn{14}{c}{dimension 7} \\
		\cmidrule{1-14}
		 vertices & 8 & 9 & 10 & 11 & 12 & 13 & 14 & 15 & 16 & 17 & 18 & 19 & 20\\
		facets min & 8 & 12 & 14 & 15 & 18 & 20 & \smash{\tikz[baseline]{\node[circle,draw, inner sep = 1pt] at (0,0.75ex) {16}}} 
																				& 39 & 55 & 67 & 100 & 139 & 219\\
		facets max & 8 & 20 & 40 & 70 & 104 & 134 & 163 & 198 & 239 & 254 & 281 & 244 & 228\\
		\cmidrule{1-14}
	\end{tabular}
	\caption{The~number of~facets of~a 2-neighborly 0/1-polytope}
	\label{tab:facets}
\end{table}

\begin{figure}[tbh]
	\centering
	\begin{tabular}{c*7r}
		 &
		\multicolumn{7}{c}{coordinates}\\
		\multirow{14}{*}{\begin{turn}{90}vertices\end{turn}} 
& 0 & 0 & 0 & 0 & 0 & 0 & 0\\
& 0 & 0 & 0 & 0 & 0 & 0 & 1\\
& 0 & 0 & 0 & 0 & 0 & 1 & 0\\
& 0 & 0 & 0 & 0 & 1 & 0 & 0\\
& 0 & 0 & 0 & 1 & 0 & 0 & 0\\
& 0 & 0 & 1 & 0 & 0 & 0 & 0\\
& 0 & 1 & 0 & 0 & 0 & 0 & 0\\
& 1 & 0 & 0 & 0 & 0 & 1 & 1\\
& 1 & 0 & 0 & 1 & 1 & 0 & 0\\
& 1 & 0 & 1 & 0 & 1 & 0 & 1\\
& 1 & 0 & 1 & 1 & 0 & 1 & 0\\
& 1 & 1 & 0 & 0 & 1 & 1 & 0\\
& 1 & 1 & 0 & 1 & 0 & 0 & 1\\
& 1 & 1 & 1 & 0 & 0 & 0 & 0\\ 
	\end{tabular}
	\caption{The~2-neighborly 0/1-polytope with 14 vertices and 16 facets}
\label{fig:p14}
\end{figure}

\section{Polytopes with a~big number of~vertices}

Let $N_2(d)$ be the~maximal number of~vertices of~a 2-neighborly $d$-dimensional 0/1-polytope.
In~\cite{Aichholzer:2000}, it was showed that $N_2(d-1) + 1 \le N_2(d) \le 2 N_2(d-1)$ and 
given some estimations for $d \le 10$.
By using Algorithm~\ref{alg:enumeration}, we improve these estimations (see Table~\ref{tab:bigd}).

\begin{table}[tbh]
	\centering
	\begin{tabular}{*{6}r}
		\toprule
		&     dimension & 7 & 8 & 9 & 10 \\
		\cmidrule{2-6}
		\multirow{2}{*}{\cite{Aichholzer:2000}}
		& best example & 18 & 25 & 33 & 44 \\
		& upper bound & 24 & 48 & 96 & 192 \\
		\cmidrule{2-6}
		\multirow{2}{*}{new}
		& best example & 20 & 28 & 38 & 52 \\
		& upper bound & 20 & 34 & 68 & 136 \\
		\bottomrule
	\end{tabular}
	\caption{The~maximal number of~vertices of~a 2-neighborly 0/1-polytope}
	\label{tab:bigd}
\end{table}



\section*{Acknowledgments}

The~author is grateful to Sergey Lichak, who helped to run Algorithm~\ref{alg:enumeration} on the~computer cluster.

%
%

\end{document}